\documentclass[12pt,leqno]{article}
\usepackage{amsmath,amssymb,amscd}
\newcommand\qed{{\unskip\nobreak\hfil\penalty50\hskip2em\vadjust{}
    \nobreak\hfil$\Box$\parfillskip=0pt\finalhyphendemerits=0\par}}

\newtheorem{theorem}{Theorem}[section]
\newtheorem{lemma}{Lemma}[section]

\newtheorem{definition}{Definition}[section]

\numberwithin{equation}{section}

\newcommand{\x}{\times}

\renewcommand{\a}{\alpha}

\renewcommand{\d}{\delta}
\newcommand{\D}{\Delta}

\newcommand{\g}{\gamma}

\renewcommand{\th}{\theta}
\renewcommand{\O}{\Omega}
\renewcommand{\o}{\omega}
\newcommand{\z}{\zeta}

\renewcommand{\i}{\infty}

\newcommand{\bN}{{\mathbb N}}

\newcommand{\bR}{{\mathbb R}}

\newcommand{\macp}{\alpha}

\begin{document}

\title{Moment Bounds for the Smoluchowski Equation and their
Consequences}
\author{Alan Hammond\thanks{This work was performed while A.H. held a postdoctoral fellowship in the Department of Mathematics at U.B.C.} \\
New York University \\
and \\
Fraydoun Rezakhanlou\thanks{This work is supported in part by NSF
grant DMS0307021}
\\
UC Berkeley}

\maketitle

\setcounter{section}{0}

\begin{abstract}
We prove $L^{\infty}\big( \mathbb{R}^d \times [0,\infty) \big)$ bounds on moments $X_a : = \sum_{m \in \mathbb{N}}{m^{a}f_m(x,t)}$ of the
Smoluchowski coagulation equations with diffusion, in any dimension $d \geq 1$. If the collision propensities $\alpha(n,m)$ of mass $n$ and mass $m$ particles grow
 more slowly than $(n+m)\big(d(n) + d(m) \big)$, and the diffusion rate $d(\cdot)$ is non-increasing and satisfies 
$m^{-b_1} \leq d(m) \leq m^{-b_2}$ for some $b_1$ and $b_2$ satisfying $0 \leq b_2 < b_1 < \infty$, then any weak 
solution satisfies $X_a \in L^{\infty}\big( \mathbb{R}^d \times [0,T] \big) \cap L^1 \big( \mathbb{R}^d \times [0,T] \big)$
 for every $a \in \mathbb{N}$ and $T \in (0,\infty)$, (provided that certain moments of the initial data are finite). As a consequence, we 
infer that these conditions are sufficient to ensure uniqueness of a weak solution and its conservation of mass.    
\end{abstract}

\begin{section}{Introduction}
The Smoluchowski coagulation equation is a coupled system of partial differential equations that describes the 
evolving densities of a system of diffusing particles that are prone to coagulate in pairs.
A sequence of functions $f_n: \mathbb{R}^d \times [0,\infty) \to [0,\infty)$, $n \in \mathbb{N}$, 
is a solution of the Smoluchowski coagulation equation if it satisfies
\begin{equation}\label{syspde}
  \frac{\partial}{\partial t} f_n \big( x,t \big) = d(n) \Delta f_n (x,t) 
  \, + \, Q_n(f)(x,t),
\end{equation}
with $Q_n(f) = Q_n^+ (f) - Q_n^- (f)$, where
\begin{displaymath}
  Q_n^+ (f) (x,t)  = \sum_{m=1}^{n-1} \macp(m,n-m) f_m(x,t) f_{n-m}(x,t) 
\end{displaymath}
and 
\begin{displaymath}
  Q_n^- (f) (x,t)  = 2 f_n(x,t) \sum_{m=1}^{\infty} \macp(n,m) f_m(x,t) .
\end{displaymath}
We will interpret this solution in weak sense. Namely, we will assume that $Q_n^+$ and $Q_n^-$ belong to
 $L^1\big(\mathbb{R}^d \times [0,T] \big)$ for each  $T \in [0,\infty)$ and $n \in \mathbb{N}$, and that
\begin{displaymath}
 f_n (x,t) = S_t^{d(n)} f_n^0 (x) \, + \, \int_0^t S_{t-s}^{d(n)}Q_n(x,s)ds,
\end{displaymath}
where $f^0$ denotes the initial data, $S_t^D$ the semigroup associated with the equation $u_t = D \Delta u$,
 and where $Q_n(x,s)$ means $Q_n(f)(x,s)$.

The system (\ref{syspde}) has two sets of parameter values, the sequence $d: \mathbb{N} \to [0,\infty)$, 
where $d(n)$ denotes the diffusion rate of the Brownian particle of mass $n$, and the collection 
$\macp: \mathbb{N}^2 \to [0,\infty)$, where $\macp(m,n)$ models the average propensity of particles of masses $m$ and $n$ to coagulate.
 The terms $Q_n^+(f)$ and $Q_n^-(f)$ are gain and loss terms for the presence of particles of mass $n$ that arise from 
the binary coagulation of particles. 

The system (\ref{syspde}) may be augmented by considering the fragmentation of particles into two or more sub-particles.
 A continuous version of the system, in which particles have real (rather than integer) mass, has also been considered.
 The data are defined $d:[0,\infty) \to [0,\infty)$ and $\alpha: [0,\infty)^2 \to [0,\infty)$, and the sums are replaced
 integrals in the definitions of $Q_n^+(f)$ and $Q_n^-(f)$. 
The spatially homogeneous version of the equations, in which each $f_n$ is a function of time alone, are better 
understood (the work \cite{ballcarr} resolved many of the central questions for the discrete case without fragmentation).

Lang and Nguyen \cite{Lang} considered a system of mass-bearing diffusing particles, whose diffusion rate was chosen to be independent 
of mass, that are prone to coagulate in pairs at close range, and demonstrated that, in a kinetic limit, the density of particles evolves 
macroscopically as a solution of (\ref{syspde}). A kinetic limit means that, with an initial number $N$ of particles, the order of
 the interaction range $\epsilon = \epsilon(N)$ of any given particle is chosen so that a typical particle experiences a rate of 
collision that is bounded away from zero and infinity for all high $N$. In \cite{hrop}, \cite{hrop2}, this kinetic limit derivation 
was extended to permit more general diffusion rates, and to include a stochastic mechanism of interaction. Suppose that particles
 of mass $n \in \mathbb{N}$ have a range of interaction given by $r(n) \epsilon$, for some increasing 
function $r: \mathbb{N} \to (0,\infty)$. That is, $\epsilon = \epsilon(N)$ determines the order of the range as a
 function of the initial particle number, whereas the function $r$ specifies the relative interaction range of particles of differing masses. 
The variations to the derivation required by introducing the radial dependence $r$ have been discussed in \cite{frpap}.
 It is shown there that, if the dimension $d \geq 3$, the macroscopic collision propensity $\macp: \mathbb{N}^2 \to (0,\infty)$
 that appears in the Smoluchowski coagulation equation satisfied by the macroscopic density profile of the particle system satisfies
\begin{equation}\label{betarel}
  \macp(n,m) \leq c \big( d(n) + d(m) \big) \Big( r(n) + r(m) \Big)^{d-2},
\end{equation}
where $c$ is a constant without dependence on the details of the stochastic interaction. 
(Note that $\beta$ rather than $\macp$ was the notation used in \cite{hrop},\cite{hrop2}. The change in notation in this paper ensures 
consistency with the PDE literature.) In fact, the left and right-hand-sides of (\ref{betarel}) are of the same order provided 
that the interaction mechanism is strong enough to ensure that a uniformly positive fraction of pairs of particles that come
 within the interaction range coagulate.
 
The growth rate of  $r: \mathbb{N} \to (0,\infty)$ presumably depends on the internal structure of the particles. 
If they are simply balls, each of the same density, then $r(n) = c_0 n^{1/d}$. On the other hand, an internal structure that 
is fractal might give rise to a relation of the form $r(n) = n^{\chi}$, for some $\chi \geq 1/d$. The physically reasonable 
range of values of $\chi$ would seem to be contained in $[1/d,1]$. This is because, whatever fractal structure a particle of
 mass $n \in \mathbb{N}$ may have, if contact with another particle is required for the pair to coagulate, then the
 interaction range $r(n)$ is at most of order $n$, (with the extremal case being that in which the particle takes
 the form of a line segment). Indeed, non-trivial fractal structure would suggest that $\chi$ is strictly less than one.

Monotonically decreasing choices of  the diffusion rates $d: \mathbb{N} \to [0,\infty)$ seem to be physically realistic,
 since the diffusive motion is presumably stimulated by the bombardment of much smaller elements of an ambient gas. 
The choice $d(n)= \frac{1}{n}$ in three dimensional space is justified if the particles are modelled as balls. 
(For a mathematical treatment,
see \cite{DGL}, which derives Brownian motion as the long-term behaviour of a ball being struck by elements in a Poisson cloud of point
 particles of Gaussian velocity.)   

In this paper, we examine the behaviour of solutions of (\ref{syspde}).
We have directed our attention to parameter values $\macp$ and $d$ that seem to be justified by
 the existing kinetic limit derivations of (\ref{syspde}) (though, in principle, other choices may arise
 from a derivation of the equations from a quite different model).

We mention firstly that, under the assumption
\begin{equation}\label{conexi}
  \lim_{m \to \infty}  \frac{\macp(n,m)}{m}  = 0  \ \ \textrm{for each $n \in \mathbb{N}$,}
\end{equation}
the global existence of a weak solution of (\ref{syspde}) has been established in \cite{LM}, which work includes 
fragmentation in the equations. From the physically reasonable assumptions that $d(\cdot)$ is uniformly bounded  and $r(n) = o(n)$, we see that 
(\ref{conexi}) is satisfied in dimension $d=3$ by choices of $\alpha$ 
satisfying (\ref{betarel}).

An important formal property of solutions of (\ref{syspde}) is mass-conservation. 
\begin{definition}
Let $f$ solve (\ref{syspde}) weakly. We say that $f$ conserves mass on the time interval $[0,T]$ provided 
that $I(t) = I(0)$ for each $t \in [0,T]$, where
\begin{displaymath}
 I(t)  = \sum_{m \in \mathbb{N}} m \int_{\mathbb{R}^d} f_m(x,t)dx.
\end{displaymath}
\end{definition}
While mass conservation holds formally, the only estimate that is readily available is $I(t) \leq I(0)$ for $t \geq 0$. Indeed,
the inequality may be strict, in which case, {\it gelation} is said to occur (at the infimum of times at which the inequality is strict).
 If this happens for a solution of (\ref{syspde}) which is obtained as a kinetic limit of a particle system, then,
 after the gelation time, a positive fraction of the mass of particles is contained in particles whose mass is greater than 
some function that grows to infinity as the initial particle number tends to infinity. 

In the spatially homogeneous setting, much progress has been made in establishing when gelation occurs. In \cite{emp} and \cite{elmp},
 continuous versions of the spatially homogeneous equations with fragmentation are considered. If the continuous analogue $\macp(x,y)$ 
of the coagulation rates $\macp$ is supposed to satisfy 
\begin{equation}\label{ctsbetarel}
\macp(x,y) = x^a y^b + x^b y^a, \, \, x,y \in \mathbb{R}
\end{equation}
with $a,b \in (0,1)$ and $a + b > 1$, then gelation occurs, unless the inhibiting effect of fragmentation is strong enough. 
It is natural to postulate from these results that if the microscopic interaction range $r(n)$ we have discussed behaves
 like $r(n) = n^{\chi}$ with $\chi> \frac{1}{d-2}$, with a diffusion rate $d$ uniformly bounded below and a bounded 
domain $\Omega$ in place of $\mathbb{R}^d$, then gelation will occur. This is because the formula (\ref{betarel}),
 (which, as already noted, may be written as an equality in the case of a reaction mechanism that is not particularly weak), 
is bounded below by the discrete analogue of the coagulation propensity given in (\ref{ctsbetarel}). As we have commented, however,
 we do not anticipate such behaviour for $\macp$ in equations arising from a three-dimensional particle system of the type 
considered in \cite{frpap}.

Rigorous sufficient conditions for mass-conservation, or for gelation, have been available in the spatially 
inhomogeneous setting only under stringent assumptions on parameters. See \cite{colpou} for the case of constant
 diffusion rates, and \cite{wrzd} for a criterion that requires uniform boundedness of $\macp$ and further 
information about the behaviour of solutions of the system (however, each of these papers includes fragmentation in the equations).
In Theorem \ref{th1.3}, we present a more applicable sufficient condition. 

The result largely depends on new moment bounds on solutions of (\ref{syspde}). Theorem \ref{th1.1} presents bounds
 on the $L^1$ norm of moments of a solution, and Theorem \ref{th1.2} provides $L^{\infty}$ bounds on such moments.   
Previously, $L^{\infty}$ estimates on solutions of (\ref{syspde}), with the effect of fragmentation included, have
 been obtained (see \cite{wrzexi},\cite{wrzd},\cite{benwrz},\cite{colpou}), under fairly restrictive assumptions on 
coefficients for coagulation and fragmentation propensity. The dependence on $i \in \mathbb{N}$ of the bounds obtained 
on $L^{\infty}$ norms of $f_i$ does not generally permit deductions about the $L^{\infty}$-norm of moments
 $\sum_{m=1}^{\infty}{m^a f_m(x,t)}$ for any $a \geq 0$ (note however that such inferences are made in \cite{colpou}
 if the diffusion rate is identically constant, or in \cite{wrzexi} if the coagulation rates $\alpha(n,m)$ decay quickly enough). 

Our final result, Theorem 1.4, provides a criterion for uniqueness of solutions of (\ref{syspde}). 
It also relies principally on the moment bounds and is 
a straightforward adaptation of the uniqueness proof of \cite{ballcarr}
which applies to the  homogeneous case.

Our results for valid for each dimension $d \geq 1$.
Each deduction, moment bound, mass conservation, or uniqueness, depends on some regularity 
in the initial data, and some assumption on the parameters of the system. We now state the various assumptions that we require.  

{\bf Assumption 1.1}
\[
\lim_{n+m\to\i} \frac {\a(n,m)}{ (n+m)(d(n)+d(m)) }=0.
\]
More precisely, for every $\d>0$, there exists $k_0=k_0(\d)>0$ such
that if $n+m>k_0$, then
\[
\a(n,m)\le \d (n+m)(d(n)+d(m)) .
\]
In addition, the function $d$ is uniformly bounded.

{\bf Assumption 1.2}\  
 The function
$d(n)$ is a non-increasing function of $n$ and that $\a(n,m)\le
C_0(n+m)$ for a constant $C_0$. Moreover, there exist positive
constants $r_1$ and $r_2$ and nonnegative constants $ b_2 \le b_1$ such that,
\[
r_1n^{-b_1} \le d(n) \le r_2n^{-b_2}.
\]

{\bf Assumption 1.3}\  The function $d$ is uniformly positive
and non-increasing, and there exists a constant $C_0$ such that
\[
\a(n,m)\le C_0(n+m) .
\]

Our notation for the various moments of $f$ will be
\begin{equation}
\label{eq1.3} X_a=X_a(x,t)=\sum_n n^a f_n(x,t),   \ \ \ \hat
X_a =\hat X_a (x,t)=\sum_n n^a d(n)^{d/2} f_n(x,t).
\end{equation}
and
\begin{eqnarray}
\label{eq1.4}
 Y_a(x,t)&=&\sum_{n,m} nm(n^a+m^a)(d(n)+d(m))
f_n(x,t)f_m(x,t),\nonumber\\
\hat Y_a(x,t)&=&\sum_{n,m} (n^{a}m+m^{a}n)\a(n,m) f_nf_m.
\end{eqnarray}
We also set
\begin{equation}
\label{eq1.5} \phi_0(x) = \begin{cases}
|x|^{2-d}&\mbox{if $d\ge 3$,} \\
-\frac {1}{2\pi} \log|x|\ 1\!\!1(|x| \le 1)&\mbox{if $d=2$,}\\
\frac 12(1-|x|)\ 1\!\!1(2|x|\le 1)&\mbox{if $d=1$.}
\end{cases}
\end{equation}

We now state the four theorems.

\begin{theorem}
\label{th1.1} {\bf {(Moment bound I)}}Assume Assumption 1.1. Then
for every $a\ge 2$ and positive $A$ and $T$,  there exists a
constant $C=C(a,A,T)$ such that, if 
\begin{equation}\label{eq1.6}
\iint \sum_{n,m} X_a(x,0) X_1(y,0) \phi_0(x-y)dxdy \le A,
\end{equation}
and
\begin{equation}\label{eq1.7}
{\mathrm ess} \sup_x\int X_a(y,0) \phi_0(x-y)dy \le A,\ \ \ \int
 X_a(x,0) dx \le A,
\end{equation}
then
\begin{equation}
\label{eq1.8} \int_0^{T} \int Y_{a-1}dx dt \le C,\ \ \ \int_0^{T}
\int \hat Y_{a-1}dx dt \le C,
\end{equation}
and
\begin{equation}\label{eq1.9}
\sup_{t\in [0,T]}\ \int X_a(x,t) dx \le C.
\end{equation}
Moreover, the constant $C$ can be chosen to be independent of $T$
when $d>2$.
\end{theorem}

\begin{theorem}
\label{th1.2} {\bf {(Moment bound II)}} Assume Assumption 1.2.
Then there exists a function $\g(a,b_1,b_2)$ with
$\lim_{a\to\i}\g(a,b_1,b_2)=\i$ such that if
 \begin{equation}\label{eq1.10}
\sum_nn^e\|f^0_n\|_{L^\i(\bR^d)}<\i,\ \ \ X_a \in L^1(\bR^d\x
[0,T]),
\end{equation}
 then
 \begin{equation}\label{eq1.11}
\sum_nn^e\|f_n\|_{L^\i(\bR^d\x [0,T])} <\i,
\end{equation}
for every $e\le \g(a,b_1,b_2)$. 
\end{theorem}

\noindent
{\bf Remark 1.1} In particular, 
if Assumptions 1.1--1.2 hold, $X_a(\cdot,0)\in
 L^1(\bR^d)$ and 
\[
\sum_nn^a\|f^0_n\|_{L^\i(\bR^d)}<\i,
\]
 for every $a\in \bN$, then $X_a \in
L^\i(\bR^d\x [0,T])\cap L^1(\bR^d\x [0,T])$  for every $a\in\bN$ and $T \in (0,\infty)$.

We refer to \eqref{eq3.15} for the explicit form of the function $\g$.
Also note that Theorem 1.1 offers sufficient conditions to ensure
$X_a \in L^1(\bR^d\x [0,T])$ (so that this theorem has been invoked in Remark 1.1).

\begin{theorem}
\label{th1.3} {\bf {(Conservation of Mass)}}Let $f$ be a weak
solution of \textup{\eqref{syspde}}.  Assume that
\eqref{eq1.8} holds for $a=2$.
Then $f$ conserves mass on the time interval $[0,T]$. 
Assume instead that $X_1(\cdot, 0)\in L^\i(\bR^d)$,
$X_2(\cdot, 0)\in L^1(\bR^d)$ and Assumption 1.3 holds. 
Then $f$ conserves mass on the time interval $[0,\infty)$.
\end{theorem}

\begin{theorem}
\label{th1.4} {\bf {(Uniqueness)}} There is a unique weak solution of
\textup{\eqref{syspde}} on the interval $[0,T]$ among those satisfying
$X_2 \in L^\i(\bR^d\times [0,T])$.
\end{theorem}

\noindent
{\bf Remark 1.2}
Assume that there exist positive constants $c_1$ and $c_2$ such that $\a(n,m)\le c_1(n^a+m^a)$ and $d(n)\ge c_2 n^{-b}$ for
all $n,m\in \bN$. Assume that $a+b<1$. As a consequence of Theorems 1.1--1.4, if $\sum_nn^e \|f^0_n\|_{L^\i(\bR^d)}<\i$
 and $\|\sum_nn^e f^0_n\|_{L^1(\bR^d)}<\i$ for 
sufficiently large $e$, then (1.1) has a unique solution which is mass conserving.

\noindent
{\bf Remark 1.3} Theorem 1.1 is also true for the continuous version of the system (\ref{syspde})
with a verbatim proof. In the continuous version, all the summations over $n$ and $m$ are replaced with
integrations with respect to $dn$ and $dm$. On the other hand, we can derive the continuous version of Theorems 1.2 and 1.3 only for a particular 
solution. The proof of Theorem 1.4 does not readily adapt to the continuous setting. See Remark 4.1 for further comments.

If the parameters $\macp$ and $d$ derived from the kinetic limit of a particle system are such that uniqueness among solutions of
 (\ref{syspde}) is unknown, then, in principle at least, the particle system may not approximate a single solution of (\ref{syspde}). 
It might, for example, approximate several different solutions, each with a positive probability. This would be very peculiar, and so,
 it is pleasing to be able to rule out the possibility by establishing uniqueness.

If we adopt the relation $r(n) = n^{\chi}$ for particle interaction range, then, recalling (\ref{betarel}), 
we have that, in $\mathbb{R}^d$ with $d \geq 3$, the macroscopic coagulation propensity arising from the microscopic random model satisfies
\begin{displaymath}
 \macp (n,m) \leq c \big( d(n) + d(m) \big) \max\big\{ n,m \big\}^{\chi (d-2)}.
\end{displaymath}
In view of Theorems \ref{th1.3} and \ref{th1.4}, we have a unique solution and mass conservation throughout 
time provided that $\chi \in \big[ 0, \frac{1}{d-2} \big)$ and suitable moments of the initial densities are finite.
 As such, the discussion following (\ref{betarel}) rules out the occurrence of gelation in three dimensions in 
 particle systems of the type considered in \cite{Lang} and \cite{hrop}.

In a similar vein to the comment that follows (\ref{ctsbetarel}), we mention that, in the case where $d$ is uniformly positive,
 Theorem \ref{th1.3} gives a sufficient condition for mass conservation of a solution of (\ref{syspde}) that is close to being sharp. 
Indeed, choices of $\macp(n,m)$ that grow much more quickly than Assumption 1.2 permits are bounded below by 
the expression in (\ref{ctsbetarel}), for some choice of $a,b \in (0,1)$ with $a + b > 1$. Such a choice of $\macp$
 would thus be expected to show gelation in the spatially homogeneous case.
Corollary 8.2 of \cite{emp} adapts the argument of the homogeneous case to assert that gelation occurs for any weak solution of the continuous 
version of (\ref{syspde}) in a bounded subset of $\mathbb{R}^d$, provided that  (\ref{ctsbetarel}) is satisfied for such $a$ and $b$ as above.  
The diffusive motion of particles 
in $\mathbb{R}^d$ may act to inhibit gelation of a solution of (\ref{syspde}), though a dense 
initial condition is likely to ensure it. 

\noindent{\bf Acknowledgment.} We thank G\'abor Pete for comments on a draft version.
\setcounter{section}{1}
\end{section}
\begin{section}{The tracer particle approach}

Each of the $L^{\infty}$ and moment bounds on solutions of (\ref{syspde}) that we present in this paper will be proved by PDE methods. 
However, for each of our results, we earlier derived a similar assertion by 
a quite different approach. A given solution of (\ref{syspde}) is understood in terms of the random trajectory of a tracer particle, 
whose behaviour is typical of the many particles that form the density profile of the solution. 
We have certainly found this random method to be intuitively appealing, and it may find application to
 other PDE for which kinetic limit derivations have been made.  We have 
thus devoted this section to explaining the tracer particle approach.

It will in fact be helpful to recall in more detail the model analysed in \cite{hrop}. A sequence of 
microscopic random models indexed by their initial number $N \in \mathbb{N}$ of particles is given. 
In the $N$-th model, each of the $N$ particles has an initial location $x(0)$ and integer mass $m(0)$ set independently according to 
\begin{equation}\label{imass}
 \mathbb{P} \big( m(0) = k \big) = \frac{\int_{\mathbb{R}^d}f_k^0(x)dx}{\sum_{n=1}^{\infty}\int_{\mathbb{R}^d}f_n^0(x)dx} ,
\end{equation} 
with $x(0)$ having law 
\begin{equation}\label{iloc}
   \frac{f_k^0 (\cdot)}{\int_{\mathbb{R}^d} f_k^0(x)dx},
\end{equation}
conditional on $m(0) = k$. At any given moment of time $t \in [0,\infty)$, particles of mass $k \in \mathbb{N}$ 
evolve as Brownian motions with diffusion rate $d(k)$, with $d:\mathbb{N} \to (0,\infty)$ a given collection of
 constants. Particles are liable to coagulate in pairs when their displacement is of order $\epsilon = \epsilon(N)$ 
(we refer the reader to the introduction of \cite{hrop} for the details of the interaction mechanism). We set 
$\epsilon \approx N^{- \frac{1}{d-2}}$ (for $d \geq 3$), to ensure that a typical particle experiences a rate
 of collision that is bounded away from zero and $\infty$ uniformly in $N$.
At the collision event, the two incoming particles disappear, to be replaced by a third, that assumes the sum
 of the masses of the ingoing two, and is located in an $\epsilon$-vicinity of either of the colliding particles.
 While the macroscopic behaviour of the system is likely to be independent of the choice of placement of the new
 particle within this microscopic vicinity of the colliding pair, it was convenient for the derivation performed 
in \cite{hrop} to assume that, if the masses of the colliding pair are $n$ and $m$, then the location of the new
 particle is taken to be that of one or other of the pair with probabilities $\frac{n}{n+m}$ and $\frac{m}{n+m}$.
 For what follows, it is convenient to regard a particle that has a collision as surviving it, and becoming the 
outgoing particle, with a probability proportional to its mass, and disappearing from the model in the other event. 

Theorem 1 of \cite{hrop} specifies the macroscopic behaviour of this particle system, when the initial particle number $N$ is taken to be high.
To summarise the result without recourse to equations, 
at typical points $(x,t) \in \mathbb{R}^d \times [0,\infty)$ of space-time, the number of particles of given
 mass $m \in \mathbb{N}$ located in a vicinity of $x$ at time $t$, normalized appropriately, approximates the
 {\em density} $f_m(x,t)$, for some solution $\big\{ f_m: m \in \mathbb{N} \big\}$
of (\ref{syspde}). 
>From the result, we may infer the law of the trajectory of a typical particle, in the limit of high initial particle number, as follows.
 
Consider a particle picked uniformly at random from the $N$ particles that are initially present. We call this the tracer 
particle in the $N$-th model.
In accordance with the preceding description, the initial mass $m(0)$ and 
location $x(0)$ of the particle are given by (\ref{imass}) and (\ref{iloc}).
We know from \cite{hrop}
that $\big\{ f_n:\mathbb{R}^d \times [0,\infty) \to [0,\infty), n \in \mathbb{N} \big\}$ gives the asymptotic particle densities 
of the microscopic models. Hence, 
asymptotically in high $N$, the trajectory of the particle is such that, if at time $t \in [0,\infty)$, the particle has 
location $x \in \mathbb{R}^d$ and mass $m$, it evolves as a Brownian motion at rate $d(m)$, and experiences collision with 
 a particle of mass $n$ at rate $\macp(m,n) f_n(x,t)$.
The details of the collision event that we specified imply that the tracer particle $(x(t),m(t))$, on colliding with a 
mass $n$ particle, survives the collision and assumes mass $m(t^+) = m(t) + n$ with probability 
$\frac{m(t)}{m(t) + n}$, and disappears from the model otherwise.

The limit in high $N$ of the law of the tracer 
particle in the $N$-th model is in fact dependent only on the given solution 
 $\big\{ f_n: n \in \mathbb{N} \big\}$
that the microscopic particle densities approximate, and not on other details of the random models, such as the location 
of other particles. We now formalise this definition: a tracer particle specified in terms of a given solution  
$\big\{ f_n: n \in \mathbb{N} \big\}$ of (\ref{syspde}) and not in terms of the data of any microscopic model.

\begin{definition}
Let a solution $\big\{ f_n: n \in \mathbb{N} \big\}$ of (\ref{syspde}) be given.
 The tracer particle governed by $f$ the random process 
$z = (x,m): [0,\infty) \to \mathbb{R}^d \times \mathbb{N} \, \cup \, \{ c \}$ whose initial law is given by 
$$
 \mathbb{P} \Big( m(0) = m \Big) = \frac{\int_{\mathbb{R}^d}f_m^0(x)dx}{\sum_{n=1}^{\infty}  \int_{\mathbb{R}^d}f_n^0(x)dx},
$$
with the initial location $x(0)$ having density $f_m^0(\cdot)/{\int_{\mathbb{R}^d}f_m^0(x)dx}$, conditional on $m(0) = m$.
At time $t$, the particle's location $x(t)$ evolves as a Brownian motion at rate $d(m)$, provided that $m(t)=m$. Moreover, 
at time $t$, and for any $n \in \mathbb{N}$, the particle is said to undertake a mass transition $m \rightarrow n + m$ at
 rate $\macp(n,m)f_n \big( x(t), t \big)$. Such a transition {\em succeeds} with probability $m/(n+m)$, in which case, 
$m(t)$ is set equal to $n+m$, and {\em fails} in the other event, in which case, the particle is relegated to the cemetery 
state $c$: $z(s)$ is set equal to $c$ for all $s \geq t$. 
\end{definition}
\noindent{\bf Remark 2.1} Certain smoothness (or measurability) assumptions are in fact required to ensure that the process 
$z = (x,m): [0,\infty) \to \mathbb{R}^d \times \mathbb{N} \, \cup \, \{ c \}$ exists, even locally in time.

Let $z: [0,\infty) \to \mathbb{R}^d \times \mathbb{N} \, \cup \, \{ c \}$ denote the tracer particle governed by a given
 solution $\big\{ f_n: n \in \mathbb{N} \big\}$ of (\ref{syspde}). Let $g_n(x,t)$ denote the density of its location at time $t$:
\begin{displaymath}
   \mathbb{P} \Big( x(t) \in A, m(t)  = n  \Big)
   = \int_A g_n \big( x,t \big) dx.
\end{displaymath}

The evolution equation for the system $\big\{ g_n : n \in \mathbb{N} \big\}$ is given by
\begin{eqnarray}
  \frac{\partial}{\partial t} g_n \big(x,t\big)
 &= &d(n) \Delta g_n (x,t) \nonumber  \\ 
& &+    \sum_{m=1}^{n-1} \macp(m,n-m) f_m(x,t) g_{n-m}(x,t) 
 \, - \, 2 g_n(x,t) \sum_{m=1}^{\infty} \macp(n,m) f_m(x,t).\label{geqn}
\end{eqnarray}
Note that the choice $g_n \equiv f_n$ solves this equation. Assuming that there is a unique non-zero solution of (\ref{geqn}), we find that
\begin{equation}\label{geqf}
  g_n \equiv f_n .
\end{equation}
Before commenting further on (\ref{geqf}), 
we want to emphasise 
how we have changed our point of view of the tracer particle. At first, we regarded it as a typical particle in a microscopic random model,
 whose behaviour in the large is described by a solution of (\ref{syspde}), and then, secondly, as a random trajectory defined purely in 
terms of such a solution. The former point of view motivates the study of the system (\ref{syspde}). We will now discuss how the latter
 is valuable in studying a given solution of (\ref{syspde}). Indeed, that the tracer particle defined by a given solution of (\ref{syspde})
 may be a useful tool for studying that solution is apparent from (\ref{geqf}): if we understand the likely behaviour of the tracer
 particle, we infer bounds on the density $\big\{ g_n: n \in \mathbb{N} \big\}$ of its location, and, by (\ref{geqf}), on the 
solution  $\big\{ f_n: n \in \mathbb{N} \big\}$ itself. 

We have stated the relation (\ref{geqf}) because doing so permits a more succinct summary of how the tracer particle approach works.
 However, in making the approach rigorous, we take a different route, which we now summarise. We construct a sequence 
$\big\{ f_n^N: \mathbb{R}^d \times [0,\infty) \to [0,\infty), n \in \mathbb{N} \big\}$, indexed by $N$, which will 
approximate the density of the tracer particle of some solution of (\ref{syspde}) when $N$ is high. 
For each $N \in \mathbb{N}$, the functions $f_n^N$, for $n \in \mathbb{N}$, are constructed inductively, on the domains 
$\mathbb{R}^d \times [0,\frac{i}{N})$, for each $i \in \mathbb{N}$. They are extended to 
$\mathbb{R}^d \times [0,\frac{i + 1}{N})$ by an inductive step in which $f_n^N(t)$ for $n \in \mathbb{N}$ and 
$t \in [\frac{i}{N},\frac{i+1}{N})$ is defined as the density of the location of a tracer particle governed by the 
constant data $f_n^N(i/N)$ during the time interval $[\frac{i}{N},t)$. For each $i \in \mathbb{N}$ and 
$T \in [0,\infty)$, the function $f_n^N$ is shown to converge in $L^1 \big( \mathbb{R}^d \times [0,T] \big)$ as $N \to \infty$
 to a limit $f_n: \mathbb{R}^d \times [0,\infty) \to [0,\infty)$. By mimicking the proof of the weak stability result of \cite{LM},
 we infer that the sequence $\big\{ f_n : n \in \mathbb{N} \big\}$ is a solution of (\ref{syspde}). Tracer particle arguments, 
of which an example will shortly be given, are then applied directly to the tracer particle densities 
$\big\{ f_n^N: n \in \mathbb{N}\big\}$, and the resulting  bounds, which hold uniformly in $N$, are 
inherited in the high $N$ limit by the solution $\big\{ f_n : n \in \mathbb{N} \big\}$.
This method for making the tracer particle approach rigorous does have the drawback of applying to only 
one solution of (\ref{syspde}). 

How may the tracer particle approach be used to prove an $L^{\infty}$-bound on the solution 
$\big\{ f_n: n \in \mathbb{N} \big\}$ of (\ref{syspde}) whose construction we have just discussed? 
To simplify the exposition, suppose that the initial condition takes the form $\int_{\mathbb{R}^d} f_1^0(x) dx = 1$,
 $f_m^0(x) = 0$ for $x \in \mathbb{R}^d$ and $m > 1$. We will now sketch a proof  of the following reformulation of 
Lemma \ref{lem3.1}. Suppose that $d: \mathbb{N} \to (0,\infty)$ is decreasing, and that $d(m) > c m^{-\frac{2}{d}(1-\alpha)}$
 for some $\alpha \in [0,1]$ and $c>0$. Then\begin{equation}\label{linfbd}
 \sum_{m=1}^{\infty} m^{\alpha} f_m(x,t) \leq C u(x,t),
\end{equation}
where $u : \mathbb{R}^d \times [0,\infty) \to [0,\infty)$ solves,
$$
 \frac{\partial u}{\partial t} = d(1) \Delta u,
$$
$u(x,0) = f_1(x,0)$. 

The solution $\big\{ f_n: n \in \mathbb{N} \big\}$ of (\ref{syspde}), being given by the density of the tracer
 particle (\ref{geqf}), has the following interpretation: for $A \subseteq \mathbb{R}^d$ a Borel set, the quantity  
$$
\sum_{m=1}^{\infty}{m^{\alpha} \int_{A} f_m(x,t) dx}
$$
is equal to the expected value of the random variable $R = R(A)$ 
equal to  $m^{\alpha}$ if the tracer particle at time $t$ has mass $m$ 
and lies inside the set $A$. 

In comparing the left- and right-hand-sides of (\ref{linfbd}), we are thus assessing the degree to which the dynamics of 
the tracer particle may increase
the expected values of the random variables $R(A)$ over those obtained by using a simple Brownian particle (for numerous 
choices of the set $A$).
For example, if the set $A$ is a small ball about $x \in \mathbb{R}^d$,
and the tracer particle is close to $x$ at some time $s$ satisfying $0 < s << t$, then a mass transition 
undertaken by the tracer particle at times shortly after $s$ 
will serve to increase the expected value of $R$, 
because the 
slower diffusion rate produced by the transition is more likely to leave the particle nearby to $x$ at the later time $t$.
 However, if a transition occurs that sharply increases the mass of the particle, then it is likely to fail, 
in which case, it contributes zero to the expected value of $R$. This latter effect limits the capacity of the tracer particle
 to focus towards $x$.  

Phrasing the question quantitatively, we ask: 
what is the expected value of $R(A)$ at time $s<t$ if a mass transition $m_1 \to m_2$ occurs at time $s$? Assuming that there 
is no other mass transition, and supposing that the tracer particle is at $y \in \mathbb{R}^d$ at time $s$, the expected value of $R$ is
$$
m_1^{\alpha} \frac{1}{\big( 2 \pi (t-s) d(m_1) \big)^{d/2}}
 \int_A \exp \Big\{ - \frac{(y-x)^2}{2(t-s)d(m_1)} \Big\} dx
$$
if no mass transition occurs, whereas, it is
\begin{equation}\label{ret}
 \frac{m_1}{m_2} 
 m_2^{\alpha} \frac{1}{\big( 2 \pi (t-s) d(m_2) \big)^{d/2}}
 \int_A \exp \Big\{ - \frac{(y-x)^2}{2(t-s)d(m_2)} \Big\} dx 
\end{equation}
if the transition does occur. (The first factor in (\ref{ret}) is the survival probability for the transition).
Given that $d(m_2) \leq d(m_1)$, the latter exponential term may be bounded pointwise by the former, and we find that
$$
\frac{\mathbb{E}\Big( R \, \textrm{if mass transition occurs} \Big)}{\mathbb{E}\Big( R \, \textrm{if it does not} \Big)}
 \leq \frac{m_1}{m_2} \frac{d(m_1)^{d/2}}{d(m_2)^{d/2}} \frac{m_2^{\alpha}}{m_1^{\alpha}}  = \frac{m_1^{1 - \alpha}}{m_2^{1-\alpha}} 
 \frac{d(m_1)^{d/2}}{d(m_2)^{d/2}}.
$$
If a whole sequence of mass transition occurs, $m_i \to m_{i+1}$ at time $t_i$ for $i \in \{1,\ldots,n \}$, with $t_i \in [0,t]$
 an increasing sequence, we similarly find that the ratio of the expected values of $R$ in the case where the sequence of mass
 transitions occurs and in that where no transition takes place, is bounded above by 
$$
 \frac{1}{m^{1-\alpha}} \frac{d(1)^{d/2}}{d(m)^{d/2}}.
$$ 
(We have omitted some details of the argument: to make it rigorous, we might use an induction on the total number of mass
 transitions undertaken prior to time $t$, and invoke the strong Markov property.)

By comparison, the expected value of $R$ if no transition occurs is $u(x,t)$, where 
$$
 \frac{\partial u}{\partial t} = d(1) \Delta u,
$$
$u(x,0) = f_1(x,0)$. By choosing $A = B(x,r)$ for each $r > 0$, we learn that
$$
\sum_{m=1}^{\infty} m^{\alpha} f_m(x,t) \leq \frac{d(1)^{d/2}}{\inf_{m \in \mathbb{N}} m^{1-\alpha} d(m)^{d/2}} u(x,t).
$$ 
>From the hypothesis $d(m) > c m^{-2/d(1-\alpha)}$, we find that
$$
\sum_{m=1}^{\infty} m^{\alpha} f_m(x,t) \leq \frac{d(1)^{d/2}}{c} u(x,t),
$$ 
as we sought.  

In fact, a more careful tracer particle argument improves this result: the same conclusion (\ref{linfbd}) may be reached under 
the weaker assumption that $d$ is decreasing and satisfies $d(m) > m^{- ( 1- \alpha ) + \epsilon}$ for some $\alpha \in (0,1)$ 
and $\epsilon \in (0,\alpha)$. This is a reformulation of our claim in Remark 1.2 which is a consequence of Theorem 1.2.

We have mentioned that each of our results has an analogue with a derivation that considers the tracer particle governed by a given 
solution of (\ref{syspde}). We will not outline the method of proof of the result on mass conservation (see Theorem \ref{th1.3}), 
stating only the characterization of gelation in terms of a tracer particle. We alter the definition of the tracer particle governed
 by $\big\{ f_n: n \in \mathbb{N} \big\}$, so that the particle survives every mass transition. Given (\ref{geqf}), mass 
conservation of the solution $f$ until a given time $T \in [0,\infty)$ occurs if and only if the tracer particle experiences 
only finitely many collisions on the time interval $[0,T]$. 
\end{section}

\section{Moment bounds under Assumption 1.1}
\label{sec2}

Let us first construct an auxiliary  function $H$ that will be
needed for the proof of Theorem \ref{th1.1}. Before doing so, let us
make an observation regarding the case $d\ge 3$. When $d\ge 3$,
define
\begin{equation}
\label{eq2.1}
H(x)=c(d)|x|^{2-d}
\end{equation}
where $c(d)=(d-2)^{-1}\o_d^{-1}$ with $\o_d$ denoting the
$(d-1)$--dimensional measure of the unit sphere
$S^{d-1}=\{x\in\bR^d:|x|=1\}$. We then have that $\Delta
H=-\delta_0$, where $\d_0$ denotes the Dirac's measure at $0$.
 More precisely, for a test
function $g$, the function
\[
u(x)=\int H(x-y)g(y)dy,
\]
satisfies $\Delta u=-g$. Note that $H\ge 0$ and this property is
lacking when $d\le 2$. Because of this we can only hope for the
existence of a suitable function $H$ such that $H\ge 0$ but now
$\Delta H=-\delta_0+Error$ for an $Error$ that can be controlled.
This is the content of our first lemma.

\setcounter{lemma}{0}
\begin{lemma}
\label{lem2.1} Assume $d\le 2$. There exist functions  $ H$ and $ K$
such that $ H \ge 0$, $ K$ is bounded, $ K$ is of compact support,
\begin{equation}
\label{eq2.2}
 -\D_x  H(x) = \d_0 -  K(x),
\end{equation}
and the function $ H- \phi_0$ is bounded. (The function $\phi_0$ was
defined in \eqref{eq1.5}.)
\end{lemma}

\bigskip
\noindent {\bf {\em Proof.}} The construction of $H$ for $d=1$ is
straightforward; we can readily find a nonnegative function $H$ such
that $H=\phi_0$ in $[-1/2,1/2]$, $H=0$ outside $[-1,1]$, and $H$ is
smooth off the origin.

For the construction of the function $ H$ when $d=2$, let us start
from the function $\phi_0$ and make an important observation. Note that
if
\begin{equation}
\label{eq2.3}
 R(x) = J * \phi_0(x) = \frac {-1}{2\pi} \int_{|x-y|
\le 1} \log|x-y|J(y)dy,
\end{equation}
with $J\ge 0$, then $R \ge 0$ and
\begin{equation}
\label{eq2.4}
-\D R = J - {\tilde J},
\end{equation}
where
\begin{equation}
\label{eq2.5}
{\tilde J}(x) = \frac {1}{2\pi} \int_{|z| = 1}
J(x-z)dS(z)=J*\tilde\delta_0,
\end{equation}
where $dS$ denotes the Lebesgue measure on the unit circle and
${\tilde \d}_0$ denotes the normalized Lebesgue measure on the unit
circle.  In \eqref{eq2.3}, we may replace $J(y)dy$ with a measure
$J(dy)$.  Then \eqref{eq2.4} is still valid weakly for an obvious
interpretation for \eqref{eq2.5}.  In particular if we choose $J(dy)
= \d_0(dy)$, then
\[
-\D \phi_0 = \d_0 - {\tilde \d}_0.
\]
Our goal is to replace ${\tilde \d}_0$ with a bounded function of
compact support.
 For this we set $\phi_1 = {\tilde \d}_0 * \phi_0$ to obtain $\phi_1
\ge 0$ and
\[
-\D(\phi_0 + \phi_1) = \d_0 - {\hat \d}_0,
\]
where ${\hat \d}_0= \tilde\d_0*\tilde\d_0$ is now a ``function'' and
is weakly given by
\begin{eqnarray*}
\int h(z){\hat \d}_0(dz)
&= &\frac {1}{(2\pi)^2} \int_{|z|=1} \int_{|a|=1} h(z+a)dS(z)dS(a) \\
&= &\frac {1}{(2\pi)^2} \int_0^{2\pi} \int_0^{2\pi} h(e^{i\th_1} +
e^{i\th_2})d\th_1d\th_2.
\end{eqnarray*}
Note that the Jacobian of the transformation $\psi:(\th_1,\th_2)
\mapsto (e^{i\th_1} + e^{i\th_2})$ is given by
$|\sin(\th_1-\th_2)|$.  Since $\cos(\th_1-\th_2) = \frac {1}{2}
|e^{i\th_1} + e^{i\th_2}|^2 - 1$, we learn that if ${\hat
\d}_0(dz)={\hat \d}_0(z)dz$ , then
\begin{eqnarray*}
\hat\d_0(z) &= & 2\frac {1}{(2\pi)^2} \frac {1}{\sqrt{1-\left( \frac
{1}{2} |z|^2-1\right)^2}} 1\!\!1(|z| \le 2) \\
&= &2\frac {1}{(2\pi)^2} \frac {2}{|z|\sqrt{4-|z|^2}} 1\!\!1(|z| \le
2).
\end{eqnarray*}
Here the factor $2$ comes from the fact that $\psi$ maps (exactly)
two points to one point because
$\psi(\theta_1,\theta_2)=\psi(\theta_2,\theta_1)$. The function
$\hat \delta$ is not bounded.  Let us apply the above procedure one
more time to define $\phi_2 = \phi_0 * {\hat \d}_0$ so that $\phi_2
\ge 0$ and
\[
-\D(\phi_0 + \phi_1 + \phi_2) = \d_0 - {\bar \d}_0,
\]
where
\begin{equation}
\label{eq2.6} {\bar \d}_0(z) = \frac {4}{(2\pi)^3} \int_{|a|=1}
\frac {1}{|z-a|\sqrt{4-|z-a|^2}} 1\!\!1(|z-a| \le 2)dS(a) .
\end{equation}
As we will see, the function $\bar \d_0$ is not bounded either.
However, the function $\bar \delta_0$ is less ``singular'' that
$\hat \d_0$. For this we show that in fact $\bar \d_0$ has a
logarithmic singularity. To do this, first observe
 that ${\bar \d}_0(z)$ is radially symmetric because both $dS$ and
$\hat \d_0$ are rotationally invariant. Hence we may assume that $z$
lies on the $x_1$-axis and $z > 0$. Observe that the integrand in
\eqref{eq2.6} is singular when either $|z-a|=0$ or $|z-a|=2$. It is
not hard to show that there exists a positive constant $c_1$ such
that if $a=e^{i\theta}$ with $\theta\in(-\pi,\pi]$, then
\begin{equation}\label{eq2.7}
c_1(|\theta| + |1-|z||) \le |z-a| \le |\theta| + |1-|z||.
\end{equation}
The first inequality will be used to treat $1/|z-a|$ singularity in
\eqref{eq2.6}. We now turn to the singularity which comes from the
factor $(2-|z-a|)^{-1/2}$. For this observe that if  $a = e^{i\th}$,
then we have
\[
\gamma(\theta):=4 - |z-a|^2 = 3 - z^2 + 2z \cos \th .
\]
Let us choose $\theta_0\in [0,\pi]$ so that
\[
3 - z^2 + 2z \cos \th_0=0\ \ \  \text{or}  \ \ \ \cos
\th_0=\frac{z^2-3}{2z}.
\]
This means that the integrand in \eqref{eq2.6} is singular at
$\th_0$ and $-\th_0$. Note that when $z\ge 0$, the singular point
$\th_0$ exists only if  $z$ belongs to the interval $[1,3]$. Also note
that if $z$ is neither close to  $1$ nor $3$, then $\th_0$ is
neither  close to $0$ nor $\pi$. From the elementary inequality
\[
\gamma(\theta)=\gamma(\th)-\gamma(\th_0)=2z(\cos
\th-\cos\th_0)=-4z\sin\frac{\th+\th_0}2 \sin\frac{\th-\th_0}2,
\]
we learn that if $z$ is neither close to  $1$ nor $3$, then we can
find a positive constant $c_2$ such that
\[
\sqrt {4-|z-a|^2}\ge c_2\sqrt{|\th-\th_0|},
\]
for $\th$ close to $\th_0$. (Recall that $z\in [1,3]$.) Since this is an integrable singularity
with respect to $d\th$-integration, we deduce that $\bar\d_0$ is
bounded if $z$ stays away from the circles $|z|=1$ and $|z|=3$.
 We now assume that $|z|$ is close to $1$. In this case $\th_0$ is
close to $\pi$ and  $|\th_0-\pi|$  is comparable to
$\sqrt{|1-|z||}$. Also, because of $|z-a|\le 2$ and
$\g(\th)-\g(\th_0)>0$ we learn that $\th<\th_0$.  We have
\begin{equation}\label{eq2.8}
{4-|z-a|^2}\ge c_3
\left[(\pi-\th)+(\pi-\th_0)\right]\left[(\pi-\th)-(\pi-\th_0)\right],
\end{equation}
for a positive constant $c_3$. When $|z|$ is close to $1$, the
integrand  in \eqref{eq2.6} is singular at $\pm \th_0$ and
``almost'' singular (see \eqref{eq2.7}) at $0$ .  From \eqref{eq2.7}
and \eqref{eq2.8} we deduce
\[
{\bar \d}_0(z) \le c_4 \left| \log\left|1-|z|\right|\right|,
\]
whenever $|z|$ is close to $1$. (Here we used the fact that
$\int_0^{\th_0} [(\pi-\th)^2-(\pi-\th_0)^2]^{-1/2}d\th$ is of order
$|\log(\pi-\th_0)|$.)

 When $|z|$ is close to $3$,
$\th_0$ is small (in fact of order $O(\sqrt{3-|z|})$), and the
condition $|z-a|\le 2$ forces $\th\in[-\th_0,\th_0]$. On the other
hand,
\[
 {4-|z-a|^2}\ge c_5|\th+\th_0| {|\th-\th_0|},
\]
for a positive constant $c_5$, implies that ${\bar \d}_0(z)$ is
bounded for $|z|$ close to $3$. Here we are using
\[
 \int_{-\th_0}^{\th_0}\left(\th_0^2-\th^2\right)^{-\frac 12}\ d\theta <\i.
\]

Putting all the pieces together, we deduce

\begin{equation}
\label{eq2.9} {\bar \d}_0(z) \le c \left[\left|
\log\left|1-|z|\right|\right| \right] 1\!\!1(|z| \le 3).
\end{equation}
 We set $\phi_3 = {\bar \d}_0 * \phi_0$ to obtain $\phi_3 \ge 0$ and
\begin{equation}
\label{eq2.10} -\D(\phi_0 + \phi_1 + \phi_2 + \phi_3) = \d_0 -  K
\end{equation}
where
\[
 K(z) = \frac {1}{2\pi} \int_{|a| = 1} {\bar \d}_0(z+a)dS(a).
\]
It is straightforward to use \eqref{eq2.9} to show  that $ K$ is
uniformly bounded. Now \eqref{eq2.2} follows from \eqref{eq2.10} by
choosing $ H = \phi_0 + \phi_1 + \phi_2 + \phi_3$. It is also
straightforward to check that the functions $\phi_1 $, $ \phi_2$,
and $ \phi_3$ are bounded. \qed

\bigskip
\noindent

Let $\zeta$ be a nonnegative smooth function of compact support with
$\int\zeta=1$ and set $\z^\d(x)=\d^{-d}\z(x/\d)$. We also define
$f^\d_n=f_n*_x\z^\d$ and $Q_n^\d=Q_n*_x\z^\d$. We certainly have
\begin{equation}\label{eq2.11}
f_n^\d(x,t)=f_n^\d(x,0)+\int_0^t d(n)\Delta
f_n^\d(x,s)ds+\int_0^tQ_n^\d(x,s)ds.
\end{equation}
Also, as it is well-known,
\begin{equation}\label{eq2.12}
\sum_n\phi(n)Q_n=\sum_{n,m}\a(n,m)(\phi(n+m)-\phi(n)-\phi(m))f_nf_m.
\end{equation}
The same identity is valid if we replace $f_nf_m$ with $(f_nf_m)*_x\z^\d$ and $Q$
with $Q^\d$. Using the fact that for $\phi(n)=n1\!\!1(n\le \ell)$,
we have $\sum_n\phi(n)Q_n\le 0$ we can readily deduce that
\begin{equation}\label{eq2.13}
 \sup_\d\sup_\ell\sup_t\int \sum_{n=1}^{\ell} nf_n^\d(x,t) dx<\i.
\end{equation}

\begin{lemma}
\label{lem2.2} Let $H$ be as in Lemma 2.1.  Then there exists a constant $c_0$ such that
\[
\sup_x \sum_{n=1}^\ell n(f_n^\d *_x H)(x,t) \le \sup_x \sum_{n=1}^\ell n (f_n^\d *_x H)(x,0) +
c_0t
\]
for every positive $\d$. We may choose $c_0=0$ when $d\ge 3$.
\end{lemma}

\bigskip
\noindent {\bf {\em Proof.}}  We have
\begin{eqnarray*}
\sum_{n=1}^{\ell} n(f_n^\d *_x H)(x,t)&=& \sum_{n=1}^{\ell} n(f_n^\d
*_x H)(x,0) -\int_0^t\sum_{n=1}^{\ell} nd(n) f^\d_n (x,s)ds\\
& &\ \ \ + \int_0^t\sum_{n=1}^{\ell} nd(n) f_n^\d *_x K(x,s)ds
+\int_0^t \sum_{n=1}^{\ell} nQ_n^\d *_x H(x,s)ds.
\end{eqnarray*}
>From the boundedness of $K$, $H\ge 0$ and $\sum_{n=1}^{\ell} nQ_n \le 0$ we
deduce
\[
\sum_{n=1}^{\ell} n(f_n^{\d} *_x H)(x,t) \le \sum_{n=1}^{\ell}
n(f_n^{\d} *_x H)(x,0) + c_1 \int_0^t \int \sum_{n=1}^{\ell} nd(n)
f_n^{\d}(x,s) dxds.
\]
We now use \eqref{eq2.13} to bound the last term
to complete the proof. \qed

\bigskip
\noindent {\bf {\em Proof of Theorem \ref{th1.1}.}}  Set
\begin{eqnarray*}
Z^\d(t)& =& \int \left( \sum_{n=1}^{\ell}
n^af_n^{\d}(x,t)\right)\left(
\sum_{n=1}^{\ell} n f_n^{\d} *_x H(x,t)\right)dx,\\
Z(t)& =& \int \left( \sum_{n=1}^{\ell} n^af_n(x,t)\right)\left(
\sum_{n=1}^{\ell} n f_n *_x H(x,t)\right)dx.
\end{eqnarray*}
We have that weakly,
\begin{eqnarray*}
\frac {d}{dt} Z^\d(t)& = &-\int
 \left( \sum_{n=1}^{\ell} n^af_n^{\d}(x,t)\right)\left(
\sum_{n=1}^{\ell} n d(n) f_n^{\d}(x,t)\right)dx \\
& &\ -\int \left( \sum_{n=1}^{\ell} n^a d(n)
f_n^{\d}(x,t)\right)\left(
\sum_{n=1}^{\ell} n f_n^{\d}(x,t)\right)dx \\
& &\ +\int \left( \sum_{n=1}^{\ell} n^a f_n^{\d}(x,t)\right)\left(
\sum_{n=1}^{\ell} n Q_n^{\d} *_x H(x,t)\right)dx \\
&&\ +  \int \left( \sum_{n=1}^{\ell} n^a Q_n^{\d}(x,t)\right)\left(
\sum_{n=1}^{\ell} n f_n^{\d} *_x H(x,t)\right)dx \\
& &\ + \int
 \left( \sum_{n=1}^{\ell} n^af_n^{\d}(x,t)\right)\left(
\sum_{n=1}^{\ell} n d(n)\ f_n^{\d}*_xK(x,t)\right)dx \\
& &+\ \int \left( \sum_{n=1}^{\ell} n^a d(n)
f_n^{\d} (x,t)\right)\left(
\sum_{n=1}^{\ell} n \ f_n^{\d} *_x K  (x,t) \right)dx \\
&=: &\O_1 + \O_2 + \O_3 + \O_4+\O_5+\O_6.
\end{eqnarray*}

We now study the various terms which appear on the right-hand side. We
certainly have
\begin{eqnarray*}
\O_1 + \O_2 &= &-\frac {1}{2} \int \sum_{n,m} 1\!\!1(n,m \le \ell)
 [n^amd(m) + n^ad(n)m + m^and(n) + m^ad(m)n]f_n^\d f_m^\d \ dx \\
&= &-\frac {1}{2} \int \sum_{n,m} 1\!\!1 (n,m \le
\ell)nm(d(n)+d(m))(n^{a-1}+m^{a-1})f_n^\d f_m^\d \ dx.
\end{eqnarray*}
>From $\sum_{n=1}^{\ell} n Q_n \le 0$, we learn that $\O_3 \le 0$. By
boundedness of $K$, \eqref{eq2.13} and the boundedness of $d(\cdot)$ we deduce that
\[
|\O_5+\O_6|\le c_1\int X^\d_a dx,
\]
where $X_a$ is defined in \eqref{eq1.3} and $X^\d_a=X_a*\z^\d$.
It remains to bound $\O_4$. Note that
\begin{eqnarray}\label{eq2.14}
\sum_{n=1}^{\ell} n^a Q_n(x,t) &= &\sum_{n,m} [(n+m)^a 1\!\!1(n+m
\le \ell) -
 n^a 1\!\!1 (n \le \ell) - m^a 1\!\!1 (m \le \ell)]\a(n,m)f_nf_m \nonumber\\
&\le &c_2 \sum_{n,m} (n^{a-1}m+m^{a-1}n)\a(n,m)1\!\!1 (n+m\le \ell)
f_nf_m =:c_2 Z_4.
\end{eqnarray}
>From this, Lemma~\ref{lem2.2}, the boundedness of $H - \phi_0$ and (\ref{eq1.7}),
\[
\O_4\le c_2\int Z_4^\d \left( \sum_{n=1}^{\ell} n f_n^{\d} *_x
H(x,t)\right)dx \le c_3 \int Z_4^\d dx,
\]
where $Z_4^\d=Z_4*\z^\d$. On the other hand, for $\d_0>0$, we can
find $k_0=k_0(\d_0)$ such that if   $k_0<\ell$, then
\begin{eqnarray}\label{eq2.15} 
 Z_4&= &\sum_{n,m} 1\!\!1 (k_0 \le n+m \le \ell) (n^{a-1}m+m^{a-1}n)
\a(n,m) f_nf_m \nonumber\\
& &+ \sum_{n,m} 1\!\!1 (n+m<k_0) (n^{a-1}m+m^{a-1}n) \a(n,m) f_nf_m \\
&\le &\d_0 \sum_{n,m} 1\!\!1 (n+m\le \ell)
(n^{a-1}m+m^{a-1}n)(n+m)(d(n)+d(m))f_nf_m\nonumber \\
&+ &2k_0^a \sum_{n,m} 1\!\!1 (n+m\le \ell)\a(n,m)f_nf_m.\nonumber
\end{eqnarray}
As a result,
\begin{eqnarray*}
\frac {d}{dt} Z^\d(t) &\le &- \frac {1}{2}  \int \sum_{n,m} 1\!\!1
(n,m \le \ell)nm(n^{a-1}+m^{a-1})(d(n)+d(m))
f_n^{\d}f_m^\d dx \\
& &+ 2 c_3 \d_0 \int \sum_{n,m} 1\!\!1 (n,m \le
\ell)nm(n^{a-1}+m^{a-1})(d(n)+d(m))
(f_nf_m)*\z^\d dx \\
& &+ 2c_3k_0^a \int \sum_{n,m} 1\!\!1 (n+m\le
\ell)\a(n,m)(f_nf_m)*\z^{\d} dx+c_1 \int X^\d_a(x,t)dx.
\end{eqnarray*}
Here, we are using the identity, valid provided that $a \ge 2$,
\begin{displaymath}
 (n+m)\big( n^{a-2} + m^{a-2} \big) \le 2 \big( n^{a-1} + m^{a-1} \big).
\end{displaymath}

We now send $\d$ to $0$ to yield
\begin{eqnarray}\label{eq2.16}
\frac {d}{dt} Z(t) &\le &\left(  2c_3\d_0-\frac {1}{2} \right) \int
\sum_{n,m} 1\!\!1 (n,m \le \ell)nm(n^{a-1}+m^{a-1})(d(n)+d(m))
f_nf_m
dx \nonumber\\
& & \ \ \ +2c_3k_0^a \int \sum_{n,m} 1\!\!1 (n+m\le
\ell)\a(n,m)f_nf_m dx+c_1 \int X_a(x,t)dx.
\end{eqnarray}
Note that the time integral of the second integral is bounded because
\[
\frac {d}{dt}\int \sum_{n=1}^{\ell} f_n\ dx\le -\int \sum_{n,m} 1\!\!1 (n+m\le
\ell)\a(n,m)f_nf_m dx,
\]
which implies,
\begin{equation}\label{eq2.17}
\int_0^T\int \sum_{n,m} 1\!\!1 (n+m\le
\ell)\a(n,m)f_nf_m dx dt\le \int \sum_{n=1}^{\ell}nf_n^0\ dx.
\end{equation}
Furthermore, the equality
\begin{equation}\label{eq2.18}
\frac {d}{dt} \int \sum_{n=1}^{\ell} n^a f_n(x,t)dx = \int
\sum_{n=1}^{\ell} n^a Q_n(x,t) dx,
 \end{equation}
and \eqref{eq2.14} imply that
\begin{eqnarray*}
\int_0^t\int \sum_{n=1}^{\ell} n^a f_n(x,s)dx d s &\le& t\int
\sum_{n=1}^{\ell} n^a f_n(x,0)dx
+c_2\int_0^t\int_0^s \int  Z_4(x,\theta) dxd\theta\\
&\le& t\int \sum_{n=1}^{\ell} n^a f_n(x,0)dx +tc_2\int_0^t\int
Z_4(x,\theta) dxd\theta.
\end{eqnarray*}
>From this, \eqref{eq2.17} and \eqref{eq2.16} we
deduce that $Z(t)-Z(0)$ is bounded above by
\[
  \left(  \big(2c_3+c_1c_2t\big) \d_0 -\frac {1}{2} \right)\int_0^t \int
\sum_{n,m\le \ell} nm(n^{a-1}+m^{a-1})(d(n)+d(m)) f_nf_m dx
ds+c_4(1+ t).
\]
We now choose $\d_0=\d_0(t)$ so that $ 1/2> (2 c_3+c_1c_2t)\d_0 $. With
this choice, the bounds in \eqref{eq1.8} follow.  From \eqref{eq1.8}, \eqref{eq2.14} 
and \eqref{eq2.18} we conclude \eqref{eq1.9}. \qed

We end this section with a variant of Theorem 1.1 that holds under
Assumption 1.3.
\begin{lemma}
\label{lem2.3}  
Under Assumption 1.3, there exists a constant $C$ such that
\begin{equation}\label{eq2.20}
\int X_2(x,t) dx \le\left( \int X_2(x,0) dx\right)\exp\left(CT\|\sum_n nf_n^0\|_{L^\infty}\ \right) .
\end{equation}
\end{lemma}
\bigskip
\noindent {\bf {\em Proof.}} We start from
\begin{eqnarray*}
\frac {d}{dt} \int \sum_{n=1}^{\ell} n^2 f_n(x,t)dx& =& \int
\sum_{n=1}^{\ell} n^2 Q_n(x,t) dx\\
&\le&\int 2\sum_{n,m} nm\a(n,m)1\!\!1 (n+m\le \ell) f_nf_m\ dx\\
&\le&\int 2C_0\sum_{n,m} nm(n+m)1\!\!1 (n+m\le \ell) f_nf_m\ dx\\
&\le&\int 4C_0\left(\sum_{n} n1\!\!1 (n\le \ell) f_n\right) \
\left(\sum_{m} m^21\!\!1 (m\le \ell) f_m\right)\ dx.
\end{eqnarray*}
This and Grownwall's inequality  imply \eqref{eq2.20}
because we can use the uniform positivity of $d(\cdot)$ and Lemma
3.1 of Section 3 to assert that $X_1 \in L^{\i}$.\qed
\bigskip

\section{Moment bounds when $d(\cdot)$ is non-increasing}
\label{sec3}

This section is devoted to the proof of Theorem \ref{th1.2}. We
start with a lemma.

\setcounter{lemma}{0}
\begin{lemma}
\label{lem3.1} Assume $d(\cdot)$ is non-increasing. Then
 \begin{equation}
\label{eqn3.1} \hat X_1(x,t)=\sum_{n=1}^{\i} n d(n)^{d/2} f_n(x,t)
\le d(1)^{d/2} u(x,t).
\end{equation}
where $u$ is the unique solution to $u_t = d(1) \D u$ subject to the
initial condition $u(x,0) = \sum_{n=1}^{\i} n f_n(x,0)$.
\end{lemma}

\bigskip
\noindent {\bf {\em Proof.}} We first establish
\setcounter{section}{4} \setcounter{equation}{1}
\begin{equation}
\label{eqn3.2} \sum_1^{\ell} n d(n)^{d/2}f_n(t) \le d(1)^{d/2}
S_t^{d(1)} \left( \sum_1^{\ell} n f_n^0\right) + d(\ell)^{d/2}
\int_0^t S_{t-s}^{d(\ell)} \left( \sum_1^{\ell} n Q_n(s)\right)ds.
\end{equation}
Here for simplicity, we do not display the dependence on the
$x$-variable. Note that $(\ref{eqn3.2})$ implies
\begin{equation}
\label{eqn3.3} \sum_1^{\ell} n d(n)^{d/2}f_n \le d(1)^{d/2}
S_t^{d(1)} \left( \sum_1^{\ell} n f_n^0\right),
\end{equation}
because $\sum_1^{\ell} n Q_n \le 0$.  Evidently $(\ref{eqn3.3})$ implies
$(\ref{eqn3.1})$.

We establish $(\ref{eqn3.2})$ by induction.  $(\ref{eqn3.2})$ is obvious when ${\ell}=
1$ by definition; in fact we have equality.
 Suppose $(\ref{eqn3.2})$ is valid.  We would like to deduce $(\ref{eqn3.2})$ with $\ell$
replaced with $\ell + 1$. To do so, first observe that if $D_1 \ge
D_2$ and $g \ge 0$, then
\begin{equation}
\label{eqn3.4} {D_1}^{d/2} S_t^{D_1} g \ge  {D_2}^{d/2}
S_t^{D_2}g.
\end{equation}
>From this and $(\ref{eqn3.2})$ we learn
\begin{equation}
\label{eqn3.5} \sum_1^{\ell} n d(n)^{d/2}f_n \le d(1)^{d/2}
S_t^{d(1)} \left( \sum_1^{\ell} n f_n^0\right) + d(\ell+1)^{d/2}
\int_0^t S_{t-s}^{d(\ell+1)} \left( \sum_1^{\ell} n Q_n(s)\right)ds
\end{equation}
because $d(\ell) \ge d(\ell+1)$ and $\sum_1^{\ell} n Q_n \le 0$.
Applying $(\ref{eqn3.4})$ to
\[
f_{\ell+1}(t) = S_t^{d(\ell+1)} f_{\ell+1}^0 + \int_0^t
S_{t-s}^{d(\ell+1)} Q_{\ell+1}(s)ds,
\]
yields
\begin{equation}
\label{eqn3.6} {f}_{\ell+1} (t) \le \left( \frac
{d(1)}{d(\ell+1)}\right)^{d/2} S_t^{d(1)} f_{\ell+1}^0 + \int_0^t
S_{t-s}^{d(\ell+1)} Q_{\ell+1}(s)ds.
\end{equation}
We multiply both sides of $(\ref{eqn3.6})$ by $(\ell+1)d(\ell+1)^{d/2}$ and
add the result to $(\ref{eqn3.5})$.  The outcome is
\[
\sum_1^{\ell+1} n d(n)^{d/2} f_n \le d(1)^{d/2} S_t^{d(1)} \left(
\sum_1^{\ell+1} n f_n^0\right) + d(\ell+1)^{d/2}
 \int_0^t S_{t-s}^{d(\ell+1)} \left( \sum_1^{\ell+1} n
Q_n(s)\right)ds.
\]
This completes the proof.\qed

\bigskip
\noindent
{\bf Remark 4.1.} For the continuous version of (\ref{syspde}), a similar proof leads to an $L^\i$ bound on
$\int_0^{\infty} nd(n)^{d/2}f_n dn$ provided that the function $d(\cdot)$ is piecewise constant, uniformly positive and 
nonincreasing. If we assume only the second and third conditions on $d(\cdot)$, we can establish the same bound on a solution provided that 
this solution can be approximated by solutions corresponding to piecewise constant $d(\cdot)$. 
Of course, if we already know the uniqueness of 
solutions to (1.1), then our $L^\i$ bound applies to all solutions. But we do not: when uniqueness is proved in Section 4, 
we will use Lemma 3.1 and its consequence Theorem 1.2.

\bigskip
\noindent {\bf {\em Proof of Theorem \ref{th1.2}.}} \\
\noindent{\bf Step 1}.
 Let us simply write $L^{\ell}$ for $L^{\ell}({\mathbb R}^d \times
[0,T])$.  We first show that if $\ell\ge 1$, then
\begin{equation}
\label{eqn3.7} {\hat X}_{a\ell+1} \in L^1 \Rightarrow {\hat X}_{a+1}
\in L^{\ell}.
\end{equation}
Indeed,
\begin{eqnarray*}
{\hat X}_{a+1} &= &\sum_n n^a\ nd(n)^{d/2}f_n = {\hat X}_1 \sum_n
n^a\
\frac {n d(n)^{d/2}f_n}{{\hat X}_1}, \\
{\hat X}_{a+1}^{\ell} &\le &{\hat X}_1^{\ell} \sum_n n^{a\ell}\
\frac {n d(n)^{d/2}f_n}{{\hat X}_1} = {\hat X}_1^{\ell-1}\ {\hat
X}_{a\ell+1},
\end{eqnarray*}
by H\"older's inequality.  This implies $(\ref{eqn3.7})$ because by Lemma
3.1, ${\hat X}_1 \in L^{\infty}$.

Recall that we are assuming
\[
r_1n^{-b_1} \le d(n) \le r_2n^{-b_2}.
\]
>From this and $(\ref{eqn3.7})$ we deduce that
\[
X_{\left(a\ell + 1 - b_2 \frac {d}{2}\right)} \in L^1 \Rightarrow
X_{\left(a+1-b_1 \frac {d}{2}\right)} \in L^{\ell}.
\]
This means that
\begin{equation}
\label{eq3.8} X_a \in L^1 \Rightarrow X_b \in L^{\frac {2a+b_2 {d}
- 2}{2b+ b_1 {d}-2 }},
\end{equation}
provided that $2b+ b_1 {d}-2 >0$ and $(2a+b_2d-2)/ (2b+b_1d-2)\ge 1$.\\
\bigskip
\noindent{\bf Step 2}.  We now try to bound $Q_n^+$ with the aid of
$\eqref{eq3.8}$.  Using $\a(n,m)\le C_0(n+m)$, we certainly have
\begin{eqnarray*}
Q_n^+ &= &\sum_{n_1+n_2=n} \alpha(n_1,n_2)f_{n_1}f_{n_2} \\
&\le &\sum_{n_1,n_2} 1\!\!1(n_1 \ge n/2 \ \mbox{ or } \ n_2 \ge
n/2)\alpha(n_1,n_2)f_{n_1}f_{n_2} \\
&\le &2C_0\left[ X_1 X_0(n/2) + X_1(n/2)X(0) \right]
\end{eqnarray*}
where $X_a(N) = \sum_{m \ge N} m^a\ f_m$.  Hence
\begin{eqnarray*}
Q_n^+ &\le &c_2n^{-\ell}[X_1 X_{\ell} + X_{1+\ell}X_0], \\
\|Q_n^+\|_{L^p} &\le &c_2 n^{-\ell}
[\|X_1 \|_{L^{\ell_1}}\|X_{\ell}\|_{L^{\ell_2}} +
\|X_{1+\ell}\|_{L^{\ell_3}} \|X_0\|_{L^{\ell_4}}]
\end{eqnarray*}
provided that $\frac {1}{p} = \frac {1}{\ell_1} + \frac {1}{\ell_2} = \frac
{1}{\ell_3} + \frac {1}{\ell_4}$. To use \eqref{eq3.8}, let us first assume that
$b_1d>2$ and that $2a+b_2d-2\ge b_1d +2\ell$.
Assume $X_a \in L^1$.  Then we use \eqref{eq3.8} to assert that
\begin{eqnarray*}
X_1 &\in & L^{\frac {2a+b_2d-2}{b_1d}},\ \ \ X_0 \in L^{\frac
{2a+b_2d-2}{ b_1d-2}} \\
X_{\ell} &\in &L^{\frac {2a+b_2d-2}{2\ell+b_1d-2}},\ \  X_{\ell+1} \in
L^{\frac {2a+b_2d-2}{2\ell+b_1d}}.
\end{eqnarray*}
Hence if we set
\begin{equation}\label{eq3.9}
 p = \frac {2a+b_2d-2}{2b_1d+2\ell-2},
\end{equation}
and assume that $p\ge 1$, $b_1d>2$, then we have that
\begin{equation}
\label{eq3.10} X_a \in L^1 \Rightarrow \|Q_n\|_{L^p} \le c\
n^{-\ell}.
\end{equation}
However, if $b_1d\le 2$, then we have that $X_0 \in L^\i$ because ${\hat X}_1 \in L^\i$.
>From this we learn that in this case \eqref{eq3.9} is true but now for
 \begin{equation}
\label{eq3.11}
p = \frac {2a+b_2d-2}{2\ell+b_1d}. 
\end{equation}
\bigskip
\noindent {\bf Step 3}.  Note that if
\[
p_D(x,t) = \begin{cases} (4\pi Dt)^{-d/2} \exp\left( -\frac
{|x|^2}{4Dt}\right) &\mbox{if $t >
0$,} \\
0 &\mbox{if $t < 0$,}
\end{cases}
\]
and $g$ is a function with $g(x,t)=0$ for $t<0$, then $\int_0^tS_{t-s}^D g(x,s)ds = (p_D * g)(x,t)$ where the
convolution is in both $x$ and $t$ variables.  Also note that
\begin{eqnarray*}
\int_0^T \int (p_D(x,t))^r dx dt &= &\int_0^T \int \left( \frac
{1}{\sqrt{Dt}}\right)^{dr} p_1\left( \frac {x}{\sqrt{Dt}},1\right)^r dx dt \\
&= & \int_0^T (Dt)^{\frac {d}{2}(1-r)} dt = c(T,r) D^{\frac
{d}{2}(1-r)}
\end{eqnarray*}
with $c(T,r) < \infty$ if and only if $r < \frac {2}{d} + 1$.  We
certainly have
\[
f_n(x,t) = (S_t^{d(n)}f_n^0)(x) + (p_{d(n)} * Q_n)(x,t).
\]
So,
\begin{equation}
\label{eq3.12} \|f_n\|_{L^{\infty}} \le \| (S_t^{d(n)}
f_n^0)(x)\|_{L^\i} + \|p_{d(n)}\|_{L^r} \|Q_n\|_{L^p}
\end{equation}
provided that $\frac {1}{r} + \frac {1}{p} = 1$.  Since $p_D \in
L_r$ with $r < \frac {2}{d} + 1$, it suffices to have
\[
\frac {1}{1+2/d} + \frac {1}{p} < 1.
\]
Choosing $p$ as in \eqref{eq3.9} or \eqref{eq3.11}  requires
\[
\frac {2b_1d + 2\ell-2}{2a+b_2d-2} \ \ {\text {or}}\ \  \frac {b_1d + 2\ell}{2a+b_2d-2} < \frac {2}{d+2}.
\]
More precisely,
\begin{equation}
\label{eq3.13} 
\ell < \begin{cases} \frac {1}{d+2} (2a+b_2d-2) - b_1d+1  &\mbox{if $b_1d>2$,}\\
\frac {1}{d+2} (2a+b_2d-2) - \frac 12 b_1d  &\mbox{if $b_1d\le 2$.}
\end{cases}
\end{equation}
In summary, we need $\ell$ to satisfy \eqref{eq3.13} and $r$ to
satisfy $\frac {d}{2}(1-r)>-1$. As a result, if $X_a \in L^1$
and $\ell$ satisfies \eqref{eq3.13}, then
\[
\|f_n\|_{L^{\infty}}\le A_n + c_3\ n^{-\ell}d(n)^{d(1-r)/2}
 \le A_n + c_4\ n^{-\ell}d(n)^{-1}
\]
where
\[
A_n = \| (S_t^{d(n)}) f_n^0(x)\|_{L^\i}\le \| f_n^0\|_{L^\i(\bR^d)}.
\]

\bigskip
\noindent {\bf Final Step}.  We have that $\sum_n n^{e}\|f_n \|_{L^{\infty}}\le \infty$ if $X_a \in L^1$,
\begin{equation}
\label{eq3.14} \sum_n n^{e-\ell} d(n)^{-1} < \infty,
\end{equation}
and,
\[
 \sum_n n^{e}\|  f_n^0\|_{L^\i(\bR^d)} < \infty.
\]
For \eqref{eq3.14} it suffices to have
\[
e - \ell + b_1 < -1.
\]
In other words,
\begin{equation}\label{eq3.15}
e \le \g(a,b_1,b_2):=\begin{cases} \frac {1}{d+2} (2a+b_2d-2) - b_1(d+1)  &\mbox{if $b_1d>2$,}\\
\frac {1}{d+2} (2a+b_2d-2) - \frac 12 b_1d -b_1-1 &\mbox{if $b_1d\le 2$.}
\end{cases}
\end{equation}
 \qed
\setcounter{section}{4}
\section{Uniqueness}
\label{sec4}

The main result of this section is Theorem~\ref{th4.1}. Theorem 1.4 is an immediate consequence of Theorem \ref{th4.1}. 

\begin{theorem}
\label{th4.1} Assume that $\a(n,m)\le c_0nm $ and let $f$ and $g$ be
two solutions with
\[
\left\| \sum_n n^2f_n\right\|_{L^{\infty}},\ \left\| \sum_n
n^2g_n\right\|_{L^{\infty} } \le A,
\]
where $L^p$ abbreviates $L^p(\bR^d\x [0,T])$. Then
\[
X(t) :=\int \sum_n n|f_n-g_n|(x,t)dx
\]
satisfies
\begin{equation}
X(t) \le e^{4c_0At}X(0),
\end{equation}
for $t\le T$. In particular, if $f_n(\cdot,0) = g_n(\cdot,0)$ for
all $n$, then $f_n(\cdot,t) = g_n(\cdot,t)$ for all $n$ and $t\in
[0,T]$.
\end{theorem}

We first state a straightforward lemma:

\begin{lemma}
\label{lem4.1} Let $u$ be a weak solution of
\[
u_t =D \Delta u + h
\]
with $u$ and $h \in L^1$.  Assume that $\psi$ is a continuously
differentiable convex function with $|\psi'(a)| \le c_1$ for a
constant $c_1$ and all $a \in {\mathbb R}$. Then
\begin{equation}
\int \psi(u(x,t))dx \le \int \psi(u(x,s))dx + \int_s^t \int
\psi'(u(x,\theta))h(x,\theta)dxd\theta
\end{equation}
whenever $0 < s < t$.
\end{lemma}

Lemma~\ref{lem4.1} is established by choosing a smooth  mollifier
$\rho_{\epsilon}$, and showing
 the inequality $(5.2)$ for $u_{\epsilon} = u * \rho_{\epsilon}$ and
$h_{\epsilon} = h * \rho_{\epsilon}$. We then pass to the limit
$\epsilon \to 0$.  We omit the details.

\bigskip
\noindent {\bf {\em Proof of Theorem 4.1.}}  Choose a continuously
differentiable convex function $\psi_{\delta}$ so that
 $\psi'_{\delta}(r) = sgn(r)$ for $r \notin (-\delta,\delta)$,
$\psi_{\delta} \in C^2$, $|\psi'_{\delta}(r)|\le 1$, $\psi_{\delta}
( 0)=0$ and $\psi_{\delta} \ge 0$.
  We then apply Lemma~\ref{lem4.1} to assert that the expression
\begin{equation}\label{eq4.3}
\int \sum_{n=1}^N n\ \psi_{\delta}(f_n(x,t)-g_n(x,t))dx,
\end{equation}
is bounded above by
\begin{eqnarray}
\label{eq4.4}
& &\int \sum_{n=1}^N n\ \psi_{\delta}(f_n(x,s)-g_n(x,s))ds
\nonumber
\\
& &\ \ + \int_s^t \int \sum_{n=1}^N
n\psi'_{\delta}(f_n(x,\theta)-g_n(x,\theta))(Q_n(f)(x,\theta)
-Q_n(g)(x,\theta))dxd\theta \nonumber \\
& &\ =\int \sum_{n=1}^N n\psi_{\delta}(f_n(x,s)-g_n(x,s))dx  \\
& &\ \ \  + \int_s^t \int  \sum_{n,m}
\alpha(n,m)(\Gamma_{n+m}-\Gamma_n-\Gamma_m)(f_nf_m-g_ng_m)dxd\theta.
\nonumber
\end{eqnarray}
where $\Gamma_n = n\ \psi'_{\d}(f_n-g_n)1\!\!1(n \le N)$.  Observe that
if $n,m \le N$, then
\[
(\Gamma_{n+m}-\Gamma_n-\Gamma_m)(f_nf_m-g_ng_m),
\]
equals
\begin{eqnarray*}
& &(\Gamma_{n+m}-\Gamma_n-\Gamma_m)(f_n-g_n)f_m +
(\Gamma_{n+m}-\Gamma_n-\Gamma_m)(f_m-g_m)g_n \\
& &\ \ \le (n+m) |f_n-g_n|f_m - n\psi'_{\d}(f_n-g_n)(f_n-g_n)f_m +
m|f_n-g_n|f_m
\\
& &\ \ \ + (n+m) |f_m-g_m|g_n + n|f_m-g_m|g_n - m\psi'_{\d}(f_m-g_m)(f_m - g_m)g_n \\
& &\le 2m|f_n-g_n|f_m +2 n 1\!\!1(|f_n-g_n| < \delta)|f_n-g_n|f_m \\
& &\ \ \ + 2n|f_m-g_m|g_n + 2 m 1\!\!1 (|f_m-g_m| < \delta)|f_m-g_m|g_n.
\end{eqnarray*}
>From this, $(5.4)$ and $\alpha(n,m) \le c_0nm$, we learn that the
expression $(5.3)$
is bounded above by
\begin{eqnarray*}
& &\int \sum_{n=1}^N n\psi_{\delta}(f_n(x,s)-g_n(x,s))dx \\
& &+ 2c_0 \int_s^t \int \left[ \sum_{n=1}^N n |f_n-g_n|\right]
\left[
\sum_{m=1}^N m^2(f_m+g_m)\right] dxd\theta \\
& &+2 c_0 \delta \int_s^t \int \left[ \sum_{n=1}^N n^2
1\!\!1(|f_n-g_n| < \delta)\right] \left[ \sum_{m=1}^N
m(f_m+g_m)\right]dxd\theta.
\end{eqnarray*}
We then send $\delta \to 0$ and $N \to \infty$, in this order to
obtain
\begin{eqnarray*}
\int \sum_{n=1}^{\infty} n|f_n(x,t)-g_n(x,t)|dx &\le &\int
\sum_{n=1}^{\infty} n |f_n(x,s)-g_n(x,s)|dx \\
& &+2 c_0 \int_s^t \int \left[ \sum_{n=1}^{\infty} n
|f_n-g_n|\right] \left[ \sum_{m=1}^{\infty} m^2(f_m+g_m)\right]
dxd\theta.
\end{eqnarray*}
The theorem now follows from this and Gronwall's inequality. \qed

\section{Mass conservation}
\label{sec5}

{\bf {\em Proof of Theorem \ref{th1.3}.}}
 We first assume that $\hat Y_1\in L^1$.
Evidently,
\begin{eqnarray*} \frac {d}{dt} \int \left( \sum_{n=1}^N n
f_n\right) dx &= &-\sum_{n,m} \{1\!\!1(n \le N < n+m)n +
 1\!\!1(m \le N < n+m)m\}\a(n,m)f_nf_m \\
&= &-2 \sum_{n,m} 1\!\!1 (n \le N < n+m)n\a(n,m)f_nf_m \\
&\ge &-2\sum_{n,m} 1\!\!1 (n \ge N/2\ \text{ or }\ m >
N/2)nm\a(n,m)f_nf_m.
\end{eqnarray*}
The limit $N\to\i$ of the time average of the right-hand side is $0$
because $\hat Y_1\in L^1$. From this we can readily deduce that
\[
\lim_{N\to\i}\left[\int \left( \sum_{n=1}^N n f_n(x,t)\right)
dx-\int \left( \sum_{n=1}^N n f_n(x,0)\right) dx\right]=0.
\]
This completes the proof when $\hat Y_1\in L^1$.

We now assume that Assumption 1.3 holds. We have,
\begin{eqnarray*}
\frac {d}{dt} \int \left( \sum_{n=1}^N n f_n\right) dx &=
&-\sum_{n,m} \{1\!\!1(n \le N < n+m)n +
 1\!\!1(m \le N < n+m)m\}\a(n,m)f_nf_m \\
&\ge &-2 C_0\sum_{n,m} 1\!\!1 (n \le N < n+m)n(n+m)f_nf_m \\
&\ge &-2C_0 \sum_{n,m} 1\!\!1 (n \le N/2,\  m > N/2) n(n+m)f_nf_m
\\
& &-2C_0 \sum_{n,m} 1\!\!1 (n > N/2)n(n+m)f_nf_m \\
&= &-\O_1 - \O_2.
\end{eqnarray*}
We certainly have,
\begin{eqnarray*}
\|\O_1\|_{L^1} &\le &c \left\| \sum_n n^2 f_n\right\|_{L^1} \left\|
\sum_{m > N/2} f_m \right\|_{L^{\i}} \\
&+  &c  \left\| \sum_n n f_n \right\|_{L^2} \left\| \sum_{m> N/2}
 m f_m \right\|_{L^{2}}, \\
\|\O_2\|_{L^1} &\le &c \left\| \sum_{n > N/2} n f_n\right\|_{L^2}
\left\| \sum_m mf_m\right\|_{L^{2}} \\
&+ &c \left\| \sum_{n > N/2} n^2 f_n\right\|_{L^1} \left\| \sum_m
 f_m\right\|_{L^{\i}},
\end{eqnarray*}
where $L^p$ abbreviates $L^p(\bR^d\x [0,T])$. Since $\sum_m mf_m \in L^2$, $\sum_m m^2f_m \in L^1$, and
$\sum_nnf_n\in L^{\i}$, by Lemmas 4.1 and 2.3, we are done.\qed

\bibliographystyle{plain}
\bibliography{Smolucho}

\end{document}